\documentclass{amsart}
\usepackage{graphicx,amssymb}
\newtheorem{thm}{Theorem}[section]

\newtheorem{prop}[thm]{Proposition}
\theoremstyle{definition}
\newtheorem{defn}[thm]{Definition}
\theoremstyle{remark}
\newtheorem{remark}[thm]{Remark}
\numberwithin{equation}{section}

\begin{document}

\title{Anchored vector bundles and algebroids}%
\author{Michel Nguiffo Boyom}%
\address{GTA  UMR CNRS 5030 Université Montpellier2}%
\email{boyom@math.univ-montp2.fr}%



\begin{abstract}Inspired by recent works of Zang Liu, Alan
Weinstein and Ping Xu, we introduce the notions of {\it CC
algebroids} and {\it non asymmetric Courant algebroids} and study
these structures. It is shown that CC algebroids of rank greater
than 3 are the same as Courant algebroids up to a constant factor,
though the definition of CC algebroids is much simpler than that
of Courant algebroids,requiring only 2 axioms instead of 5.  The
situation is similar to that of Lie algebroids, where in the usual
definition used by all of he experts there
is a {\it redundant axiom}, e.g.[GG,KO1,KO2,MK,PL].\\
Non asymmetric Courant algebroids are shown to be nothing but
(pseudo)clan bundles (in the sense of E.B. Vinberg-Katz) which
arise in affine geometry of convex bounded domains. The study of
CC algebroids and non asymmetric Courant algebroids involves the
cohomology theory of Koszul-Vinberg algebras and their modules.

\end{abstract}
\maketitle
\section{INTRODUCTION}
 Let $M$ be a connected smooth manifold and $V$ be a real vector
 bundle on $M$. The main subject of the present paper is
 the study of smooth vector bundles with extra algebraic or geometric structures . \\
 Let $\Gamma (V)$
 be the real vector space of smooth sections of $V$. A vector bundle morphism
 $\rho$
 from $V$ to the tangent bundle of the base manifold $M$ is called an
 anchor. The vector
 space $\Gamma(V)$ is a (left) module of the associative commutative algebra
 $F(M) : = C^{\infty}(M) $
 of smooth real valued functions on $M$, but in general the multiplication
 defining
 a real algebra structure of $\Gamma(V)$ is not required to be $F(M)$-
 bilinear. The role of the anchor map is to control relationships between
 the algebra structure of $\Gamma(V)$  and its $F(M)$-module structure.
 That is the main idea behind geometric objets such as Lie algebroids,
 Koszul-Vinberg
 algebroids, Lie-Rinehart algebras and Courant algebroids. \\
 On the other hand the
 anchor $\rho$ induces a linear map from the vector space $\Gamma(V)$ to
 the vector space  $X(M)$ of smooth vector fields on the base manifold
 $M$. \\
 There are two situations depending on whether the multiplication map in
  $\Gamma(V)$ is skew
  symmetric or not.\\
   To each skew symmetric multiplication, (that we denote by
   $[.,.]$),
  one assigns the so
  called Jacobi anomaly, namely

   $$ J(s,s',s") = \oint [[s,s'],s"]$$

  where $\oint $ denotes the cyclic sum in $ s, s', s" $ \\

  In the present paper, multiplications which are not skew symmetric will be
  called
  $\textit{non }$  $\textit{asymmetric}$. To each non asymmetric multiplication
  we will assign its {\it  Koszul-Vinberg anomaly}, namely

   $$KV(s,s',s") = (s,s',s")-(s',s,s")$$

 where  $ (s,s',s") =  s.(s'.s")-(s.s').s" $ is the associator
 (which vanishes for associative algebras). \\

 From a non asymmetric multiplication, say $s.s'$, one
 can construct a skew symmetric one by setting

   $$ [s,s'] = s.s'-s'.s $$

 The Jacobi anomaly of the last bracket is related to the
 Koszul-Vinberg anomaly  by the following equation

   $$ J(s,s',s") = \oint K(s,s',s") $$

 Given an element $s$ of $\Gamma(V)$, its image $\rho(s)$ under the anchor
 map acts on $F(M)$ as first order
 differential operator.
 The relationship between the real algebra
 structure of $\Gamma(V)$ and its $F(M)$-module structure is
 controlled by the following {\it Leibniz anomaly}

    $$ \textit{L}(s,f,s') = s.(fs') - (\rho(s)f)s'-f(s.s') $$

 We will call an $\textit{almost}$  $\textit{algebroid}$ on the base manifold $M$
 any couple  $(V,.)$  consisting of a vector bundle $V$
 on $M$ together with a real algebra structure  $(\Gamma(V),.)$
 in the vector space of smooth sections of $V$.\\

 The present work is concerned with the study of the couple
 $ (J(s,s',s"),L(s,f,s'))$ (resp. $(KV(s,s',s"),L(s,f,s'))$)
 of Jacobi anomaly and Leiniz anomaly (resp. the
 Koszul-Vinberg anomaly and Leibniz anomaly)
 of
 an anchored almost algebroid whose multplication is skew symmetric (resp. non
 asymmetric). \\
 For instance one easily sees that

 $$ (J(s,s',s"),L(s,f,s')) = (0,0) \,\, \forall s, s' \in
 \Gamma(V)\, \forall f \in F(M) $$

 if and only if $ (V,\rho,[.,.]) $ is a Lie algebroid.
 On the other hand

 $$( KV(s,s',s"),L(s,f,s') = (0,0) \,\, \forall s, s', s" \in
 \Gamma(V),\, \forall f\in F(M) $$

 if and only if $(V,\rho,.)$ is a Koszul-Vinberg algebroid.\\

 Our work is inspired by those of Z.J. Liu,
 A. Weinstein and P. Xu on Dirac structures, (see [LWX1,LWX2].)\\
  We have adopted
 the use of the cohomology of the algebra $F(M)$ viewed as a Koszul-Vinberg algebra .\\
  Resolving a old problem raised by
 Gerstenhaber [GM], we recently constructed  the cohomolgy theory of
 Koszul-Vinberg algebras and their modules which controls
 deformations of those structures, [NB1]. That cohomology is also related to
 Poisson geometry [NB2,NB3]. \\
 It is
 remarkable that
 from each vector bundle $V$ on the base manifold $M$ arise two cochain complexes
 of the Koszul-Vinberg algebra $F(M)$, namely

  $$\leqno(c)\quad C^{\ast}(F(M),V ) =
  \bigoplus_{k}Hom(\bigotimes^{k}F(M),\Gamma(V)). $$

  $$\leqno(c^{\star}) \quad  C^{\ast}(F(M),V^{\ast}) =
  \bigoplus_{k}Hom(\bigotimes^{k}
    F(M),\Gamma(V^{\star})). $$

  To control the Jacobi anomaly and the Leibniz anomaly, the cohomology
 theory of the Koszul-Vinberg $F(M)$ turns out to be more efficient than the
 Hoschschild cohomology of the associative algebra $F(M)$. \\
 It is reasonable to conjecture
 that many ingredients that are involved in the theory of
 Courant algebroid structures  and Dirac structures lie in the derived objets of the complexes (c) and
  $(c^{\star})$. \\
 In the present work, we will introduce the notions of $\textit{CC
 algebroids}$ and $\textit{non asymmetric}$
 Courant algebroids. The cochain complexes (c) and $(c^{\star})$
 will be used to study these structures. \\

 Our main result concerning CC algebroids is Theorem 5.1
 which in particular implies that the system of $\textrm{five}$
 $\textit{axioms}$ in the usual definition of $\textit{Courant}$
 $\textit{algebroids}$ contains three axioms which are superfluous whenever the
 rank of the vector bundle is greater than 3. For more details
 on Courant
 algebroid structures and related topics, the reader may consult works of
 Liu-Weinstein-Xu, mainly [LWX1,LWX2].\\

 We take this opportunity to recall that some years ago (1995 and 2000)
 we have pointed out a similar redundancy in the usual
 definition of Lie algebroid structures. Namely:  a {\it Lie algebroid }is an
 anchored almost Lie algebroid $(V,\rho,[.,.])$ such
 that the following axioms hold
 $$\leqno(AX1) \quad J(s,s',s") = 0. $$
 $$\leqno(AX2) \quad L(s,f,s') = 0. $$
 $$\leqno(AX3)  \quad \rho([s,s']) = [\rho(s),\rho(s')] = 0. $$

  However (AX3) is {\it superfluous }.\\
  In fact, it is easily seen that (Ax3) is a consequence of the other two axioms
  (AX1) and (AX2).
 The reader is referred to Section 2 of our paper [NB1] in  Banach Center Publications ,
 Vol. 54, page 103, Warszawa 2001; ibidem, page 45, DEFINITION 2, joint paper by J.Grabowski and K.
 Grabowska, [GG],
 contains the
 superfluous axiom (AX3). Thereafter, the later joint paper by J.Grabowski and
 G. Marmo, [GMa], attests that those authors hadn't read the reference we just
 recalled. \\

  Another consequence of our Theorem 5.1 is that the theory of
  Courant
  algebroid structures of rank less than three
 differs from that of Courant algebroid structures of rank greater than three.
 That
 phenomenon is illustrated by our Example 5.2\\

 Non asymmetric Courant algebroid structures are studied in Section 7 and Section 8.
 They lead us to
 features which are quite different from those inherited from CC algebroid
 structures. Indeed, those structures lead to locally flat (pseudo)clan bundles
 (Theorem 7.2).
 That phenomenon is an unexpected incursion of the
 affine geometry of homogeneous bounded domains in the theory of non
 asymmetric Courant algebroid structures. We
 obtain interesting  relationships between the differential geometry of non
 asymmetric
 Courant algebroid structures and their Koszul-Vinberg cohomology (Theorem 7.5).\\
  An another result concerns the class of non asymmetric Courant
  algebroid structures with a definite forms $<.,.>$. We will show that under
  some
  additional conditions, the cohomology class of $<.,.>$
  doesn't vanish (Theorem 8.1).
 This fact is in contrast to properties of clans which arise from the affine
 geometry of homogeneous hyperbolic bounded  domains,
 [KV,JLK1,KJ].
 (The cohomology theory used in [KJ] is derived from the Chevalley-Eilenberg
 cohomology theory of Lie algebras). Nevertheless Theorem 8.1  may be compared
 to similar results in [JLK2] where Jean-Louis Koszul has pointed out a lot of
 canonical vector bundle valued superorder differential forms whose cohomology
 classes never vanish. (For instance, the divergence class associated to a volume form ,
 cocycles defined by torsion free linear connection).
 Two examples of clan bundle and pseudo clan bundle are given in
 Section 8. Section 9 is devoted to some miscellaneous items. In section 10 some observations are made about relationships of algebroid structures
  with various topics.

 \section{ALGEBROIDS}

  Given a connected smooth manifold $M$, the associative algebra of smooth real
  valued functions defined on $M$ is denoted by $F(M)$.\\

  Let $V$ be a smooth vector bundle on $M$ and let
  $\Gamma(V)$ be the vector space of smooth sections of $V.$ We shall
  consider $\Gamma(V)$ as a Koszul-Vinberg module of $F(M)$ by
  setting the following axioms

      $$ \leqno(1) \qquad \qquad (sf)(x) = (fs)(x) = f(x)s(x)  $$

 for any $s$ in $\Gamma(V)$ and any  $f$ in $F(M).$
 A vector bundle $V$ will be called an {\it almost algebroid}
 whenever $\Gamma(V)$ is
 endowed with a real algebra structure. Therefore the product of two sections $s$ and $s'$
 will be denoted by $ss'$. \\

 Given an almost algebroid, it is to be noticed that in general $\Gamma(V)$ is not
 an algebra over the ring $F(M).$ \\
 A vector bundle $V$ on $M$ together with a vector bundle morphism $\rho$ to
 the tangent bundle $TM$ is called anchored vector bundle. The anchor $\rho$ induces a
 map from  $\Gamma(V)$ to $\Gamma(TM)$ which is $F(M)$ linear.
 The anchor map of an anchored almost algebroid $V$ is used to relate the
 $F(M)$-module structure of $\Gamma(V)$ with its real algebra structure. Such relationships
 yield to the concept of labelled algebroids.\\
 Before pursuing, let us recall some important labelled almost
 algebroid structures .\\

 (e1): A Lie algebroid is an anchored almost Lie algebroid
$(V,\rho,[.,.])$
 such
 that

 $$ \leqno(2i)\qquad (\Gamma(V),[.,.])\hbox{ is a real Lie algebra. }$$

     Given $s$, $s'$ in $\Gamma(V)$ and $f$ in $F(M)$ one
 has

 $$ \leqno(2ii)\qquad [s,fs'] = (\rho(s)f)s' + f[s,s']. $$

 (e2): A {\it Koszul-Vinberg algebroid }is an anchored almost algebroid
 $(V,\rho,.)$ such that given elements $s$, $s'$ and $s"$
 of $\Gamma(V)$ and an element $f$ of $F(M),$ one has

  $$ \leqno(3_{i})\quad   s.(s'.s") - (s.s').s" - s'.(s.s") + (s'.s).s" = 0. $$

      $$ \leqno(3_{ii}) \quad (fs).s' - f(s.s') = 0. $$
     $$  \leqno(3_{iii}) \quad s.(fs') - (\rho(s)f)s' - f(s.s') = 0.$$

 \section{ALMOST LIE ALGEBROIDS}

 In this section, we will be concerned with the so called {\it almost Lie
 algebroid } structures, viz those almost algebroids $(V,[.,.])$ whose
multiplications
 $[.,.]$ are skew symmetric.
 Let $(V,[.,.])$ be an anchored almost Lie algebroid on $M$. Let $s$, $s'$ and $s"$ be
 sections of $V$ and let $f$ be an element of $F(M)$. The only obstructions for
 an almost Lie algebroid $(V,\rho,[.,.])$ to be an Lie algebroid
 are\\

 (ob1): Jacobi anomaly

 $$ J(s,s',s") = \oint[[s,s'],s"]. $$

 (ob2): Leibniz anomaly

   $$ L(s,f,s') = [s,fs']-(\rho(s)f)s'-f[s,s']. $$

 Regarding the case of non asymmetric anchored almost algebroid
 structure, we will replace the Jacobi anomaly by the following
 quantity,which is called {\it Koszul-Vinberg anomaly}:

$$ \leqno(4) \quad  KV(s,s',s") = s.(s'.s") - (s.s').s" - s'.(s.s") +
(s'.s).s". $$

 We intend to point out that the cohomology theory of Koszul-Vinberg algebras
 and their modules provides tools which are useful in studying  the
 Jacobi anomaly and the Koszul-Vinberg anomaly.
 This idea has been inspired to us by the theory of Courant algebroid
 structures. [LWX1,LWX2,LW], see also [UK].\\

 \section{ KV-COHOMOLOGY  $H^{\ast}(F(M),V)$}

 Recall that an algebra $A$ whose associator is symmetric with
 respect to
 the first two
 arguments, viz $ KV(a,b,c) = 0 $   $ \forall a, b, c\in A $, is
 called a {\it Koszul-Vinberg
 algebra}. In particular any associative algebra is a Koszul-Vinberg algebra
 . So
 is the case for
 $F(M)$ when it is endowed with its natural associative commutative
 real algebra structure.\\
 A two-sided module
 of $F(M)$, say $W$, is called a {\it Koszul-Vinberg module} if
 the following identities hold

 $$f(gw)-(fg)w = g(fw)-(gf)w $$
 $$f(wg)-(fw)g = w(fg)-(wf)g,\, \forall f, g \in F(M),\,  \forall w\in W $$

 Let $V$ be a vector bundle on $M$. Then, according to (1), the
 vector space
 $\Gamma(V)$ is a Koszul-Vinberg module of $F(M)$.
 We shall deal with the cochain complex whose $k^{th}$
 homogeneous space is the vector space  $C^{k}(F(M),V)$ of $k$-multi-linear maps from
 $F(M)$ to $\Gamma(V)$, $k$ being a positive integer.
 When $k$ = 0
 we set

         $$C^{0}(F(M),V) = \Gamma(V)$$

 The coboundary operator

           $$\delta : C^{k}(F(M),V)\rightarrow
C^{k+1)}(F(M),V)$$

is defined as follows

 $$\leqno(5_{i}) \quad    \delta = 0 \quad if \,\,k = 0. $$

 If $k$ is a positive integer, then

 $$\leqno(5_{ii}) \quad \delta(\Theta)(a_{1},..,a_{k+1}) =
\sum_{j}(-1)^{j}((a_{j}\Theta)(a_{1},..,\hat{a}_{j},..,a_{k+1})+$$
  $$ a_{k+1}(\Theta(a_{1},..,\hat{a}_{j},...,a_{k},a_{j}))) $$

  where

  $$ (a_{j}\Theta)(a_{1},.,a_{k}) = a_{j}(\Theta(a_{1},...,a_{k}))
  -\sum_{r}\Theta(a_{1},...,a_{j}a_{r},...,a_{k})$$

 It is easy to check that

  $$ H^{0}(F(M),V)=\Gamma(V), $$

  $$ H^{1}(F(M,V)=Der(F(M,\Gamma(V)).$$

 Thus, two cocycles in $C^{1}(F(M),V)$ are cohomologuous if and only if there are
 equal.\\

 Let $(V,\rho,[.,.])$ be an anchored almost Lie algebroid and let us suppose
 that the
 corresponding vector bundle is endowed with a symmetric bilinear
 form which is denoted by $<.,.>$. We adopt notations of [LWX1,LWX2].  To
 each triple ($s$,$s'$,$s"$) of elements of $\Gamma(V)$ we assign the smooth
 function
 $T(s,s',s")$ which is defined by

   $$T(s,s',s") = \oint<[s,s'],s">. $$

 The high lighted focus in the theory of Courant algebroid structures consists
 of using the function $T$ to control the Jacobi anomaly (see [LWX1,LWX2]).\\

 \section{THE MAIN THEOREM}

 Keeping in mind the notations used above, we are in position to prove
 the following statement.\\

 \begin{thm}
 \label{theorem:main} Let $(V,[.,.],\rho )$ be an anchored almost Lie algebroid
 on $M$ .
 Let one suppose that the following assumptions to hold.\\
 (i) The vector bundle $V$ is endowed with a non degenerate
 symmetric bilinear form which is denoted by $<.,.>$.\\
 (ii) There is a cocycle $D$ in $C^{1}(F(M),V)$ satisfying the following
 two identities
 $$ \leqno(r1) \quad J(s,s',s")=D(T(s,s',s")) \quad \forall s, s', s" \in
 \Gamma(V). $$
 $$\leqno(r2)  \quad \rho(s)(<s',s">) = <[s,s']+D(<s,s'>),s"> +
  <s',[s,s"]+D(<s,s">)>. $$
 If $rank(V) > 3,$ then the anchor map $ \rho $ satisfies the following
 identity\\
 $$ \rho([s,s'])=[\rho(s),\rho(s')]  $$
 \end{thm}

 {\it Proof}\\ Following [UK], the hypothesis (r2) allows one to control the
Leibniz
  anomaly. More precisely let $f \in F(M)$. Taking in account both (r2)
  and the the $\delta$-closeness of $D$, a direct calculation of the quantity
  $\rho(s)(<fs',s">)$
  yields to the following identity

$$ \leqno(6)     \quad [s,fs']-(\rho(s)f)s'-f[s,s'] = -<s,s'>D(f).$$

Let $s$ and $s'$ be fixed elements of $\Gamma(V)$. Under the
assumption that
 $rank(V) > 3,$ we can choose a non zero element $s"$ of $\Gamma(V)$ such
 that

 $$ <s,s"> = <s',s"> = <[s,s'],s"> = 0. $$

 Therefore, for each $f \in F(M)$  the Jacobi anomaly $J(,s,s',fs" )$ is reduced to

$$ \leqno (7) \quad J(s,s',fs") = fJ(s,s',s")+T(s,s',s')D(f)+
(\rho([s,s'])f-[\rho(s),\rho(s')]f)s". $$

 On the other hand, under the same hypothesis as above, one easily checks
 the following identity

 $$T(s,s',fs" ) = fT(s,s',s"). $$

 By the virtu of (r1), the following identity holds

 $$ J(s,s',fs" ) = D(T(s,s',fs")). $$

 Combining those results with the closeness assumption $\delta(D) = o,$  we
 conclude that the following quantity

 $$ \rho([s,s'])f-\rho(s)(\rho(s')f)+\rho(s')(\rho(s)f)$$

 vanishes identically. That ends the proof of Theorem 5.1 $\square$\\

  {\it EXAMPLE 5.2}\\

 Theorem 5.1 fails when the $rank(V) < 3.$ \\
 Indeed, let $\mathbf{M}$ be the field  of real numbers. Let us set

 $$V = M \times R. $$

 Elements of $V$ are denoted by $(x,y_{x})$
 where $x$ and $y_{x}$ are two real numbers.\\
 On the other hand, let us denote the tangent bundle of $M$ by

 $$ TM = M \times R\partial_{x}.$$

 Smooth sections of $V$ are real valued smooth
 functions of one real variable. Let $f$, $g$ and $h$ be three real
 valued smooth functions
 defined on $M$. Let us define the bilinear symmetric form on $V$ by
 setting

  $$<f,g>(x)=f(x)g(x)\,\, \forall f, g \in \Gamma(V).$$

 We define the almost Lie algebroid structure on $V$ by the following bracket

  $$ [f,g]=f\partial_{x}g -g\partial_{x}f\,\, \forall f, g \in \Gamma(V). $$

 We now define the anchor map $\rho$ on $\Gamma(V)$ by putting

  $$\rho(f)=2f\partial_{x}. $$

 The 1-cocycle $D \in C^{1}(F(M),V)$ is defined by

  $$D(f)=\partial_{x}f. $$

 The reader will easily verify that the data just defined, say
 $(V,\rho,[.,.],<.,.>,D),$ satisfy both conditions
 (r1) and (r2) of
  Theorem 5.1. Nevertheless it is easily seen that the conclusion of
  Theorem 5.1 fails $\lozenge$\\

 \section{CC ALGEBROIDS }

 Considerations to be discussed in this section are inspired by Theorem 5.1 and some problems
 which are raised in [LWX2] and in [UK]. \\

\begin{defn}
\label{defn:CC}
 A {\it CC algebroid }is a datum $(V,\rho,[.,.],<.,.>,D)$
 where $(V,\rho,[.,.],<.,.>)$ is an anchored almost Lie
algebroid endowed with a non degenerate bilinear symmetric form
$<.,.>$ and $D$ is 1-cocycle in $C^{1}(F(M),V)$ with relationships
(r1) and (r2) stated in Theorem 5.1, namely \\
$$ \leqno(r1) \quad J(s,s',s") = D(T(s,s',s")\,\, \forall s, s', s" \in \Gamma(V).$$
$$ \leqno (r2) \quad \rho(s)<s',s" > = <[s,s']+D(<s,s'>),s">+<s',[s,s"]+D(<s,s">)> $$
\end{defn}

The notion of CC algebroid structure that we just introduced is
different from that of Courant algebroid structures studied  by
Lu, Weinstein and Xu.(See [LWX1],LWX2]). Below is the usually
given definiton of Courant algebroid structures \\

\begin{defn}([LWX1,LWX2]) A {\it Courant algebroid }is an anchored almost
Lie algebroid, say $(V,\rho,[.,.]$, endowed with a non degenerate
symmetric 2-form $<.,.>$ and with a 1-cocycle $D \in
C^{1}(F(M),V)$ subject to satisfy the following five
 axioms\\
 $\forall s, s', s"\in \Gamma(V), \forall f \in F(M)$ the following identities hold \\
 $$\leqno Ax1 \quad 3J(s,s',s")=D(T(s,s',s").$$
 $$ \leqno Ax2 \quad  \rho([s,s'])=[\rho(s),\rho(s')].$$
 $$\leqno Ax3 \quad [s,fs']=(\rho(s)f)s'+f[s,s']-<s,s'>D(f).$$
$$ \leqno Ax4 \quad  \rho(D(f))= 0. $$
$$\leqno Ax5 \quad  \rho(s)<s',s"> =
<[s,s']+D(<s,s'>),s">+<s',[s,s"]+D(<s,s">)> $$
\end{defn}

\begin{remark} Our Theorem 5.1 implies that up to a constant
factor, each
 CC algebroid of $rank >3$ is a Courant algebroid. Our assertion is made clear
 by the following Proposition which is a straight corollary of Theorem
 5.1.\end{remark}

\begin{prop}
 Let $(V,\rho,[.,.],D,<.,.>)$ be a CC
 algebroid whose rank is greater than three. Then, $\forall s, s', s" \in \Gamma(V)$ and
 $\forall f\in F(M)$  the following identities hold
 $$\leqno(i) \quad [s,fs'] = (\rho(s)f)s'+f[s,s'] - <s,s'>D(f)$$
 $$\leqno(ii) \quad \rho([s,s']) = [\rho(s),\rho(s')]$$
 $$\leqno(iii) \quad \rho(D(f)) = 0. $$ \end{prop}

 {\it Proof}\\ By the virtue  of (r2), (Theorem 5.1), a direct calculation
 of $\rho(s)<fs',s"> $
yields to Identity (i). Identity (ii) is nothing but the
conclusion of Theorem 5.1. To end the proof,
 one only calculates the following expression

   $$ \rho([s,fs']) = (\rho(s)f)\rho(s') + f\rho([s,s']) -
   <s,s'>\rho(D(f))$$

 Taking into account that Identity (ii) holds, one easily checks the following

   $$<s,s'>\rho(D(f)) = 0 $$

Proposition 6.4 is proved.  $\square$ \\

 Here is an another direct consequence of Definition 6.1 :

\begin{thm}
 Given a CC algebroid $(V,\rho,[.,.],<.,.>,D)$,
 the following assertions are equivalent: \\
 $$ \leqno(A1)  \quad \rho([s,s'])=[\rho(s),\rho(s')]\,\, \forall s, s' \in
 \Gamma(V)$$
 $$ \leqno(A2) \quad \rho(D(f))=0\,\, \forall \ f \in F(M).
 $$\end{thm}

 {\it Proof.}\\
 First. (A2) implies (A1) \\

 Step1. If $rank(V) = 1,$ then, let us choose $s \in \Gamma(V)$
such that $s$ is a basis of the $F(M)$-module $\Gamma(V)$ in an
open subset $U\subset M.$ Therefore, let $V_{U}$ be the inverse
image of $U$ under the projection of $V$ on $M.$ Then, $\forall
s'\in \Gamma(V_{U}),$ $\exists f\in F(M)$ such that $s'=fs.$ The
Leibniz equation gives the following identity $$
[s,s']=[s,sf]=(\rho(s)f)s-<s,s>D(f).$$ On the other hand, we have
$$ [\rho(s),\rho(fs)]=[\rho(s),f\rho(s)]=(\rho(s)f)\rho(s)$$ By
the virtue of (A2) we have $$\rho([s,fs]=(\rho(s)f)\rho(s).$$
In conclusion, (A2) implies (A1) if $rank(V)=1.$\\

 Step2. Suppose that $rank(V)>1$. Let $s, s', s" \in \Gamma(V)$
 and $f\in F(M).$Our hypothesis is that $ \forall f\in F(M)$ one
 has $\rho(D(f)) = 0.$ .
 Then, the calculation of $\rho(J(s,s',fs")) $ yields to

  $$ \rho(J(s,s',fs")=(\rho([s,s'])-[\rho(s,\rho(s')])f)\rho(s")+$$
 $$ <s',s">\rho([s,D(f)]) -<s,s">\rho([s',D(f)]). $$

  Since the left
 member of the equality above vanishes, we deduce the following identity

   $$\leqno(\ast) \quad (\rho([s,s']) -[\rho(s,\rho(s')])f)\rho(s") =
 <s,s">\rho([s',D(f)]) -<s',s">\rho([s,D(f)]). $$

 Now, let us choose an element $g\in F(M)$ satisfying the following
 two conditions in some open sub-set of the base manifold $M$

  $$\leqno (C1)  \quad <s,D(g) = 0. $$

  $$\leqno(C2)  \quad  <s',D(g)> \neq 0. $$

 Therefore, replacing $s"$ by $D(g)$ in $(\ast)$ we obtain the
 following identy\\

   $$\leqno(\ast\ast) \quad <s',D(g)>\rho([s,D(f)] = o \,\, \forall f\in F(M) $$

 Thus, the right member of the identity $(\ast)$ vanishes identically. \\

 Second: (A1) implies (A2) \\

 Now, our assumption is that $\rho$ is an algebra homomorphism from $(\Gamma(V),
 [.,.])$ to the Lie algebra of smooth vector fields on the base manifold $M$.
 Then, from the following {\it Leibniz equation}

 $$ L(s,f,s') = - <s,s'>D(f),$$

 one easily deduces that  $\rho(D(f))$ vanishes identically. That ends the proof of
 Theorem 6.5   $\square$ \\

 N.B. In [UK], Uchino raises the question to know whether the axiom (Ax2) of
 Courant algebroid structures may be deduced from the other axioms. Example 5.2
 and Theorem 6.5
 show that this question is a relevant one. Theorems 5.1, Proposition 6.4 and
 6.5 give the complete answer to $Uchino^{,s}$ question. \\
 However the two axioms (Ax3) and (Ax4) in the usual definition of Courant algebroid
 structures are {\it always superfluous.} \\
 On the other hand the three
 axioms {\it (A2),(A3) and (A4) are superfluous whenever the
 rank
 of the Courant algebroid is greater than three.} \\
 The author recently brought Alan $Weinstein^{,s}$ attention to the last observations .
 Our Theorem 5.1 shows that only the two axioms (Ax1) and (Ax5) are necessary to define Courant algebroid structures  of $rank >
 3$.\\

 So, in regard to a Courant algebroid structure, say
 $(V,\rho,[.,.],<.,.>,D),$ the cases  where $rank(V) \leq 3 $
 are quite different from those where $rank (V) > 3)$. \\
 In the cases where $rank(V) < 3$, it becomes necessary to
 add the axiom (Ax2) ( or its equivalent $ \rho(D(f)) = 0 \,
   \forall f\in F(M). $)

 Many years ago(in 1995 and in 2000) we pointed out a similar remark about
 the system of three axioms in the usual definition Lie algebroid structures.
 The correct definition of Lie algebroid structures is that we have written
 out , [NB2]. Let us  recall it below.)\\

 \begin{defn} A {\it Lie algebroid }on the base manifold $M$ is an anchored almost Lie algebroid
 $(V,\rho,[.,.])$ on $M$ with the following two properties  \\

 $$\leqno(P1) \quad  J(s,s',s" ) = o \,\forall s,s', s" \in \Gamma(V). $$

  $$\leqno(P2) \quad  [s,fs']-(\rho(s)f)s'-f[s,s'] = 0\, \forall f \in
  F(M).$$
  \end{defn}

\begin{remark}

 Both properties (P1) and (P2) imply that the anchor map $\rho$ induces a Lie algebra
 homomorphism from $(\Gamma(V),[.,.])$ to the Lie algebra of smooth vector
 fields on
 the base manifold.\end{remark}

 Regarding the {\it abundance }of literature on the theory of Lie algebroid
 structures
 we concluded and claimed (in 1995) that the redundancy of the axiom (Ax3),
 namely \\

            $$ \rho([s,s']) = [\rho(s),\rho(s')] $$

 has remained {\it unknown }to the totaly of experts for many decades. Today in
 our knowledge the contrary is still {\it uncertain}. That is reason why, once more, we would
 like to repeat things here. First authors to be recently
 convinced are J.P. Dufour, A. Banyaga, J. Leslie, T.Z. Nguyen,
 A.Weinstein  [private communications]; J. Grabowski and M. Marmo, [GMa].\\

{\it Digressions}.

 Regarding various generalizations of the theory of Lie Algebroid
 structures, the only exiting problem is to
 handle the Lie algebroid structure defect. That defect is represented by
 the couple consisting of Jacobi anomaly and Leibniz anomaly of
 anchored almost Lie algebroid structures. That is the main concern of many
 fundamental works. For instance [KO1,KO2,LWX1,LWX2,LX,MK]. To handle
 the Lie algebroid structure defects,
 many interesting ideas arise from [PP]. \\

 The highlighted point behind the theory of Courant algebroid structures is to
 ask both Jacobi anomaly and Leibniz anomaly to lie
 in the kernel of the anchor map,(via some special first order differential
 operator $D$, which is really a 1-cocycle of the complex (5) (of the Koszul-Vinberg
 algebra $F(M)$)).\\
 Similar ideas work in anchored almost Koszul-Vinberg algebroid structures.
 In the next section we intend to perform the idea that Courant
 algebroid structures provide an
 efficient framework for many interesting investigations, (see[LWX1,LWX2] for more
 details about other relationships, (such as Manin triple, Dirac structures
 and so on).\\

 \section{ NON ASYMMETRIC COURANT ALGEBROIDS}

  We plan pointing out close relationships between non
  asymmetric almost algebroid structures, (viz those $(V,.)$ such that
  the multiplication of the real algebra  $(\Gamma(V),.)$ is not assumed to be
  skew symmetric) and  the geometry of some class of bounded domains.\\

 Let $(V,\rho,.)$ be an anchored almost algebroid on the smooth manifold $M.$
 To elements $s$, $s'$ and $s"$ of $\Gamma(V)$ is assigned the associator
 $s(s's"-(ss')s"$ where $ss'$ stands for $ s.s'$. Let us recall that
 $(V,\rho,.)$ is a Koszul-Vinberg algebroid if the  following two axioms hold

  $$ \leqno(kv1) \quad s(s's") - (ss')s" - s'(ss")+(s's)s" = 0\,\,\, \forall
  s, s', s" \in \Gamma(V).$$

 $$\leqno (kv2) \quad s(fs')-(\rho(s)f)s -fss' = 0\,\, \forall f \in
  F(M). $$

  We recall that the Koszul-Vinberg anomaly is the following quantity \\

         $$(s,s',s") - (s',s,s") $$

  where $(s,s',s")$ stands for $ s(s's")-(ss')s". $\\

 Given an anchored almost algebroid on the base manifold $M$, say
 $(V,\rho,.)$,
 its KV-algebroid structure defect is represented  by the couple consisting of
 the Koszul-Vinberg anomaly and the following Leibniz anomaly

  $$  \textit{L}(s,f,s') = s(fs')-(\rho(s)f)s'-f(ss'),$$

 KV stands for Koszul-Vinberg. \\
 Let us consider an anchored almost algebroid with a non degenerate symmetric
 bilinear form,
 say $(V,\rho,.,<.,.>) $ \\

 We shall consider the vector space $\Gamma(V)$ endowed with its $F(M)$-module
 structure defined by (1).
 Now, let us set the following definition \\

 \begin{defn} A {\it non asymmetric Courant algebroid }is an anchored almost
 algebroid $(V,\rho,.)$ with a non degenerate symmetric bilinear form, say
 $<.,.>$,
 and  with a 1-cocycle $D \in C^{1}(F(M),V)$ subject to satisfy the following
requirements: \\ $\forall s, s', s" \in \Gamma(V),$ \,\,  $\forall
f \in F(M)$ one has

 $$\leqno(R1) \quad (s,s',s")-(s',s,s") = D(\delta(<.,.>(s,s',s")).$$

 $$\leqno (R2) \quad (fs)s' = f(ss').$$

 $$\leqno (R3) \quad  \rho(s)<s',s">  =  <ss'+D(<s,s'),s"> + <s',ss"+D(<s,s">)>, $$

 the  right member of the first equality in (R1) has the following meaning

        $$ \delta<.,.>(s,s',s") = -\rho(s)<s',s"> + <ss',s"> + <s',ss"> + $$
        $$ \rho(s')<s,s"> -<s's,s"> - <s,s's"> $$\end{defn}

 Our first result concerning  non asymmetric Courant algebroid structures is the following
 statement. \\

 \begin{thm} Let $(V,\rho,.,<.,.>,D)$ be a non asymmetric Courant algebroid.
 If its rank is greater than two, then the anchor map $\rho$ satisfies the following
 identity

     $$[\rho(s),\rho(s')] = \rho(ss')-\rho(s's) \,\, \forall s,
     s'\in  \Gamma(V)$$ \end{thm}

{\it Proof}\\
Let $s$ and  $s'$  be elements of $\Gamma(V)$ and let $f$ be an
element of $F(M)$. Then the following identity is a straight
consequence of (R3)
\\

      $$s(fs') - (\rho(s)f)s' - f(ss') = - <s,s'>D(f). $$

 Thus, (R3) is an efficient tool to handle
 the Leibniz
 anomaly. Since the rank of $V$ is greater than two let $s"$ be a non zero
 element
 of $\Gamma(V)$ such that

 $$<s,s">  =  <s',s"> = 0.$$

 Therefore, using the identity we just
 pointed out, a direct calculation yields to the following identity

     $$ (s,s',fs") - (s',s , fs" ) = f((s,s',s") -(s',s,s")) +
 (\delta<.,.>(s,s',s">))D(f) + $$
 $$(([\rho(s),\rho(s')] - \rho(ss') + \rho(s's))f)s". $$

 On the other hand, a  similar calculation yields to the following identity

  $$\delta<.,.>(s,s',fs") = f\delta<.,.>(s,s',s"). $$

 Therefore, by the virtu of (R1) one must conclude that the
 quantity

  $$([\rho(s),\rho(s')] - \rho(ss') + \rho(s's))f $$

 vanishes identically. That ends the demonstration of
 Theorem 7.2 $\square$\\

 Let us make some remark. Let $V,.,<.,.>$ be an non asymmetric Courant
 algebroid and let $s, s', s"\in \Gamma(V)$. As above, let us
 put

  $$KV(s,s',s")=(s,s',s")-(s',s,s"). $$

 Then $V,.,<.,.>$ gives rise to the anchored
 almost Lie algebroid structure $V,[.,.]$ whose bracket is
 defined by

 $$[s,s']=ss'-s's. $$

 The Jacobi anomaly of the last almost Lie algebroid structure is related to the
 Koszul-Vinberg anomaly $KV(s,s',s")$
 as follows

   $$ J(s,s',s") = \oint KV(s,s',s"). $$

 The digressions above lead to close relationships between non asymmetric
 Courant algebroid structures on a base manifold $M$
 and locally hessian Lie group bundles on the same base
 manifold $M$. \\
 In fact, consider a non asymmetric Courant algebroid
 $(V,\rho,.,<.,.>,D)$. Let us use (R1), (R2) and (R3) to calculate
 the quantity $\rho(fs)<s',s">$. Then, we obtain the following identity

  $$<s,s'><D(f),s"> + <s,s"><D(f),s'> = 0. $$

  Therefore, we must conclude that  $D = 0 $. The last condition is equivalent to
 $\rho = 0 $. Thus, a non asymmetric Courant algebroid is
 nothing but a Koszul-Vinberg algebra bundle endowed with a non degenerate symmetric
 bilinear form
 which is invariant under the left multiplication by elements of
 $\Gamma(V)$.
 We can write out those particular items in terms of the real valued
 cohomology the complex (5).\\

 Roughly speaking, let $\textbf{R}$ be an associative commutative
ring and let $A$ be a
 $\textbf{R}$-Koszul-Vinberg algebra. We will endow $\textbf{R}$ with the trivial
 $A$-module structure. We now consider the cochain complex whose the $k^{th}$ homogeneous
 subspaces is the vector space

 $$ C^{k}(A,K) = Hom_{K}(\otimes^{k}A,K).$$

  The coboundary operator is defined as in (5). \\

  Considering the case of a non asymmetric Courant algebroid on
  the base manifold $M$, say
  $ (V,.,<.,.>)$, we are dealing with a cohomology class in
  $H^{2}(\Gamma(V),F(M))$ containing a non degenerate cocyle, namely $<.,.>$.\\

 Let us return to the general case of non asymmetric Courant algebroid
 structures. Let
 $(V,.,<.,.>)$ be
 such an algebroid structure. We consider elements $s,s',s" \in \Gamma(V)$
 and an element $f \in F(M)$. By the virtue of (R3) one has the following identity\\

 $$  <ss,s" > + <s',ss" > = 0.  $$

 Thus, regarding the bilinear form $<.,.>$ as an element of $C^{2}(\Gamma(V),F(M)),$
 one easily sees that

      $$ \delta<.,.>(s,s',s") = 0. $$
 An interesting consequence of the last calculations is the
 following statement.\\

 \begin{thm} Each non asymmetric Courant algebroid on the base manifold $M$
 is a locally flat (pseudo) clan bundle on $M.$ \end{thm}

 {\it Proof}\\
 Let us recall that (by definition) a real clan (resp. pseudo clan) is a
 couple $(A,<.,.>)$ of real Koszul-Vinberg algebra $A$ together with a
 positive definite (resp non degenerate) real valued symmetric 2-cocyle
 $<.,.> \in C^{2}(A,R ),$ [VK,VE,SH]. \\
 A (pseudo) clan $A,<.,.>$ is
 {\it locally flat} when the left multiplication by each element of
 $A$ lies in the orthogonal algebra of $<.,.>$. \\
 Considering the case of a non asymmetric Courant algebroid  the
 vanishing property of
 the anchor map implies that the associator map $(s,s',s")$ is symmetric
 with respect to the pair ($s$,$s'$). Therefore, we get the following identity

 $$ KV(s,s',s") = o\,\, \forall s, s', s" \in \Gamma(V). $$

  Moreover, if  $x$ is a fixed element of the base manifold $M$,
  then  $\forall s, s' \in \Gamma(V)$, the element $(ss')(x)$ of $V_{x}$
  depends on $s(x)$ and on $s'(x)$
  only. Thus if we set

  $$ s(x).s'(x) = (ss')(x) $$

 then, the fiber $V_{x}$ is a Koszul-Vinberg algebra endowed with a non degenerate
 symmetric 2-cocycle, namely $<.,.>(x) $. That ends the proof of Theorem 7.3 $\square$\\

 \begin{remark} Let us keep in mind the conclusion of  Theorem 7.3, the  question
 rises to know whether the Koszul-Vinberg
 algebra bundle deduced from a non asymmetric Courant algebroid is locally trivial.
 In other words is there a {\it Koszul-Vinberg algebra} fiber type for the bundle
 $(V,.)$? \end{remark}
From the theoretic viewpoint, the cohomology theory of
Koszul-Vinberg
 algebras is helpful in studying this question. To perform
 the last idea, one must remind that the cochain complex to be considered is
$C^*(\Gamma(V),V)$ whose  coboundary
  operator is recalled below. Let $\Theta$ be an element of $C^{k}(\Gamma (V),V)$
  and let $s_{1},..,s_{k+1}$ be smooth sections of $V$, then

 $$\leqno (8) \quad \delta\Theta(s_{1},.,s_{k+1}) =
 \sum_{j}(-1)^{j}((s_{j}\Theta)(s_{1},.,s_{j-1},s_{j+1},.,s_{k+1})+$$
 $$(\Theta(s_{1},.,s_{k},x_{j}))s_{k+1})$$

    Following  our previous remarks, (see the demonstration of Theorem 7.2,)
    the coboundary operator $\delta$ is $F(M)$-linear. Thereafter a helpful tool
    in  answering the question raised in Remark 7.4 lies in  $ H^{2}(V_{o},V_{o}),$
    where $V_{o}$ stands for a fixed fiber of the vector bundle $V.$ \\
 In fact, the deformation theory of Koszul-Vinberg algebras may be controlled
 by cohomology classes of the complex (8). So, our
 digressions allow the application of a classical {\it rigidity
 theorem},[KM,GM1,KM]. More precisely, we can state the following result\\

 \begin{thm} Let $(V,.,<.,.>)$ be a non asymmetric Courant algebroid on
 a connected base manifold $M$. If $H^{2}(V_{x}, V_{x})$ vanishes
 $\forall x \in M$, then
 the K-V algebra bundle $(V,.)$ is a locally trivial.\end{thm}

{\it An outline of Proof}\\
 Without loss of generality, we may suppose the vector bundle $V$ to
 be a trivial bundle.
 Since $M$ is connected, given arbitrary points $ x_{o}, x \in M$, there is an
 isotopy  $(V_{x(t)},.)$ whose extremities are $(V_{o},.)$ and
 $(V_{x},.)$; $V_{o}$ stands for the fiber of $V$ at the point $ x_{o}$. Under the vanishing hypothesis,
 i.e. $H^{2}(\Gamma(V),V)=0 \,\, \forall x\in M,$ all of
  the fibers $(V_{x},.)$ is
 isomorphic to the fixed Koszul-Vinberg
 algebra
 $(V_{o},.)$. Let us denote by $KV(V_{o})$ the set of Koszul-Vinberg
 algebra structures on the vector space  $V_{o}$. We denote by $ \mu_{o}$ the
 Koszul-Vinberg multiplication
 that $V_{o}$  inherits from $(V,.,<.,.>)$. Under the action in
 $ Hom(\otimes^{2}V_{o},V_{o}) $ of the
 linear group of the vector space $V_{o}$, the orbit of $\mu_{o}$ is a Zariski
 open subset of $KV(V_{o})$.
 Those ingredients are used to obtain smooth family $\phi_{x}$ of isomorphisms
 from
 $(V_{x},.,)$ to $( V_{o},\mu_{o})$ $\square$\\

  \section{ A NON VANISHING THEOREM }

  Let $(V,\rho,.,<.,.>,D)$ be a non asymmetric Courant algebroid on the base manifold
  $M$. According to Theorem 7.3, such a datum may be regarded as
  a (pseudo) clan bundle on $M. $ \\
  Keeping notations in Section 7, we denote by $G_{x}$ the
  connected and simply connected Lie Group whose Lie algebra is the vector
  space $V_{x}$ endowed with the bracket defined by
  $$ [s(x),s'(x)]= (ss')(x) - (s's)(x).  $$
  Under some additional conditions, a relevant non trivial invariant of
  $(V,.,<.,.>$ is the cohomology class of the bilinear form $<.,.>.$ \\
   To make precise our assertion, let us set the following definition \\

   \begin{defn} A non asymmetric Courant algebroid is called
   {\it co-compact} if each Lie group $G_{x}$ contains a
   co-compact lattice, say $\Lambda_{x}$ \end{defn}

   Many homogeneous convex domains are base manifolds of co-compact non
   asymmetric Courant algebroids, [KJL3,KV,VEB].\\
   Below, we are going to perform that idea.\\

 \begin{thm} Let $(V,.,<.,.>)$ be a co-compact non asymmetric
 Courant
 algebroid. If the cocycle $<.,.> $ is definite, then its cohomology
 class in $H^{2}(\Gamma(V),F(M))$ doesn't vanish.\end{thm}

 {\it Proof}\\
 First of all, if the multiplication in $(\Gamma(V),.)$
 is the zero map, then the
 conclusion of Theorem 8.2 holds. Now let us suppose that the multiplication
 in $(\Gamma(V),.)$ is not the
 zero map. Let us assume the cocycle $ <.,.>$ to be (positive definite and)
 exact. Then there is a
 1-cochain $\Theta \in C^{1}(\Gamma(V),F(M))$ such that

 $$<.,.> = \delta\Theta. $$

 In other words, one has

   $$ <s,s'> = \Theta(ss') \, \forall s, s'\in  \Gamma(V).$$
  We know that $<.,.>$ is invariant
 under the left multiplications by
 elements of $\Gamma(V)$. Let $s$, $s'$ and $s"$ be elements
 of $\Gamma(V)$. For each $ x \in M,$ let
 $ G_{x}$ be the connected and simply connected Lie group whose Lie algebra
 is the vector space $ V_{x}$ endowed with the bracket
 defined by

   $$ [s(x),s'(s)] = (ss')(x) - (s's)(x). $$

 Let $\nabla$ be the left invariant linear connection on $G_{x}$
 defined by

 $$(\nabla_{s}s')(x) = (ss')(x). $$
 Actually, the differential form $\Theta $
 is De Rham closed. Then the locally flat manifold $(\textit{G}_{x},\nabla )$
 carries the
 (left) invariant closed 1-form $\Theta_{x}$ whose covariant derivative,
 say $ \nabla (\Theta)$, is positive definite. Both $\Theta $ and $\nabla $ are
 left invariant in each Lie group $G_{x}$. Therefore, the triple
 $( \Theta,\nabla,\Lambda_{x})$ gives rise to a {\it hyperbolic locally flat
  structure }on the manifold $\Lambda_{x}\backslash G_{x} $
 whose simply
 connected covering is the triple $ (G_{x}, \Theta,\nabla).$
 Therefore, each $(G_{x},\nabla )$ is isomorphic to a convex cone
 not containing any straight line. Thereafter, following [KJL1], the manifold
 $G_x$
 carries a (unique) smooth vector field $H$ satisfying the following identity \\

  given any smooth vector field $X \in \Gamma(TG_{x})$ the following identity hold
  $$\nabla_{X}(H)=X. $$

 Thereafter, let us consider elements $s,s' \in V_{x}$ as left
 invariant vector fields on $ G_{x}$. Since $\nabla $ is the Levi-Civita
 connection of the Riemannian structure $(G_{x},<.,.>)$, we check that the
 following identity holds

 $$ <ss',H > + <s',\nabla _{s}(H) > = 0. $$

Thus, the vector field $H$ is left invariant in the Lie group
$G_{x}$. From the last identity we deduce that

$$<ss',H> = - <s',s>.$$

From the exactness of the 2-cocycle $<.,.>$, we deduce the
following identities

   $$ <ss',H > = \Theta ((ss')H) = \Theta(ss'). $$

$$ <s',\nabla _{s}(H)> = <s', sH> = <s',s> = \Theta (s's). $$

In conclusion, we deduce from the calculations above the following
identity

 $$ \Theta (ss') = <s,s'> = 0\,  \forall s,s' \in V_{x}. $$

That is absurd and ends the proof of Theorem 8.2 $\square$\\

  Example  8.3 \\

  Let $L$ be the linear endomorphism of $R^{2}$
  defined by

 $$ L(x,y)=(y,x)\,\, \forall (x,y)\in R^{2}. $$

 We now consider the semi-direct product of $R^{2}$ with the one
 parameter subgroup generated by $L$. We obtain the connected and simply connected
 Lie group whose Lie algebra
  is $R^{3}$ endowed with
 the following bracket

 $$[(x,y,z),(x',y',z')] = (zy'-z'y,zx'-z'x,o). $$

 That Lie group carries a left invariant locally flat structure
 defined by the following left invariant linear connection

 $$  \nabla _{(x,y,z)}(x',y',z') = (zy',zx',o). $$

Actually, each pair $(\alpha,\beta )$ of real numbers with
$\alpha\beta\neq 0$ defines the following left invariant metric

  $$ <(x,y,z),(x',y',z')> = \alpha(xx'-yy')+\beta(zz'). $$

It is easily seen that the metric defined above is a non exact
cocycle ${\lozenge}$\\.

Example 8.4\\

Let us consider Lie algebra structure in $R^{3}$ defined by the
following bracket

  $$[(x,y,z),(x',y',z')] = (zy'-z'y,zx'-z'x,0). $$

 The associated connected and simply connected Lie group ,say $G$ , carries a
 left invariant locally flat structure corresponding to the
 following multiplication

 $$(x,y,z).(x',y',z') = (zy',-zx,0).$$

 If $\alpha$ is a non zero real number, then we define the
 following on exact cocycle

   $$ <(x,y,z),(x',y',z')> = xx'+yy'+\alpha (zz').$$

 The connected and simply connected Lie group associated to the
 Lie algebra which is defined above contains $Z^{3}$ as a co-compact lattice ${\lozenge}$ . \\

 \section{MISCELLENEA}

 Let $(V,.,<.,.>)$ be a non asymmetric Courant algebroid on $M.$ Once for
 all, let us fix an element $x_{o}$ of the base manifold $M.$
 We regard $(V,.,<.,.>)$ as a smooth deformation of the
 (pseudo) clan $(V_{o},<.,.>)$, where
 $V_{o}$ stands for the fiber of $V$ at $x_{o}.$
 As in Section 8, to each $(V_{x},<.,.>)$ is assigned the
 connected and simply connected Lie group $G_{x}$  whose Lie
 algebra of is
 the vector space $V_{x}$ endowed with the bracket defined by

$$  \leqno(9)  \quad [s(x),s'(x)]=(ss')(x)-(s's)(x).$$

 All of those Lie groups is endowed with a left invariant locally flat
 structure, (equivalently, with a left invariant locally flat linear
 connection, say $\nabla.$) Each $G_{x}$ also
 carries a left invariant locally hessian (pseudo) Riemannian  metric. That
 picture forms a smooth deformation
 of $(G_{o},\nabla,<.,.>),$ where $\nabla_{x}$ is the covariant
 derivation in $G_{x}$
 defined by the Koszul-Vinberg multiplication in $V_{x}.$ So, we can
 view the multiplication $\mu_{x}$ in
 each Koszul-Vinberg algebra {$(V_{x},.)$ as new multiplication on the
 same fixed
 vector space $V_{o}.$ Therefore the cochain $\nu_{x}$ = $\mu_{x}-\mu_{o}$ is
 a
 Koszul-Vinberg element of the complex $C^{\ast}(V_{o},V_{o}). $ In other
 words, $\nu_{x}$ satisfies the following KV equation

$$ \leqno(10) \quad  \delta(\nu_{x}) + KV_{\nu_{x}} = 0. $$

 The complex in consideration is that in REMARK 7.4.
 Koszul-Vinberg elements of that complex are the analogues of the classical
 Maurer-Cartan
 elements  which
 arise from the deformation theory of associative algebra structures  and Lie
 algebras structures. The equation (10) above is the analogue of the Maurer-Cartan equation \\

      $$ \delta\nu+1/2[\nu,\nu] = 0. $$

  (See [GM, NA, NR, KM, LWX1, [VI] and other references ibidem).\\

 We recall that to each $\nu \in C^{2}(V_{o},V_{o})$ is assigned the cochain
  $ KV_{\nu}\in C^{3}(V_{o},V_{o})$ which defined
 by

    $$ KV_{\nu}(s,s',s") = \nu(s,\nu(s',s"))-\nu(\nu(s,s'),s")-\nu(s',\nu(s,s")) +
  \nu(\nu(s's),s"). $$

  To end those miscellaneous items, let us denote by $\textit{G}$ the
  union of all of the $G_{x}$ when $x$ runs over the base
  manifold $M$. Then, $\textit{G}$ is a set bundle over
  $M$ under the set projection

  $$ G_{x} \longrightarrow x. $$

   We equip $\textit{G}$ with
  the finest topology that makes open the projection we just
  defined. We observe that $(G_{x},\nabla_{x},<.,.>_{x})$ depends smoothly
  on $x$. So, we obtain the locally hessian Lie group bundle $(\textit{G},\nabla,<.,.>)$
  on the base
  manifold $M$.\\
  Naturally, arises the question to know whether $(\textit{G},\nabla, <.,.> )$ is
  a locally trivial bundle. The complex (8) is an ingredient for
  studying the affinely flat Lie group bundle $(\textit{G},\nabla,).$
  In particular, under the hypothesis of Theorem 7.4,
  $(\textit{G},\nabla)$ is a locally trivial {\it affine Lie group} bundle on
  the base manifold $M.$ \\

  \section{OBSERVATIONS}
  (O1) Clans arose from the geometry of convex
  domains, [KV],KJL3]. In particular, the Lie algebra of a
  locally simply
  transitive
  group of affine transformations of a convex cone containing no straight line is
  a clan.
  More details can be found in fundamental papers by E.B. Vinberg,
  e.g. [EBV]. See also [KJL1,KJL2,KJL3,SH,VJ].\\

  (O2) The literature on the theory of Lie algebroid structures is impressive.
  We have related some aspects of that theory with the cohomology
  theory of Koszul-Vinberg algebras. In regard to global invariants of
  algebroid structures
   those relationships are efficient ([NB1,NB2,NBW1,NBW2]).\\

  There are many other aspects, such as the third  Lie Theorem,[AM,DP].
  The holonomy and the monodromy principle, the duality theory are studied.
  The
  theory of Singularities, and so on. There is an
abundance
  of references, for instance [DC,BR,HJ,DJ,DP,DV,DZ,NTZ MK PJ,PL,WA,WeA]. \\

  (O3) Relationships with Poisson structures and singular foliations are
  exciting
  also and have been widely studied from various viewpoints. For
  instance [FR,Ia,KJ] deal with characteristic classes viewpoint.
  The normal forms are the aim of [DJ,DZ,NBW2,NTZ]. Under some careful
  subtle techniques, the last viewpoint also walks in the theory of KV-algebroid
  structures, [NBW2].\\

  (O4) Above, we just mentioned that the theory of characteristic
  classes of Lie algebroid structures is subject of
  intense
  research programs,[FR,KJ]. Given a CC algebroid on the base manifold $M,$
  say
  $(V,\rho,[.,.],<.,.>),$ let
  $\textit{F}$ be the image of $\Gamma(V)$ under the anchor map
  $\rho.$ If the rank of $V$ is greater than three, then
  $\textit{F}$ is a subalgebra of the Lie algebra $X(M)$ of smooth
  vectors fields on $M.$ Unfortunately the Frobenius theorem generally fails
  for singular differential systems, [AM]. However it may occur that
  $\textit{F}$ be completely integrable
  in the sense of Stefan; in such an occurrence the techniques developed by
  R.Fernandes, J.Kubarski (among others) will provide characteristic
  classes of
  $(V,\rho,[.,.]),$ though the last triple fails to be an Lie
  algebroid. See [FR] for a similar remark on Courant algebroid structures.
  For instance let $\Pi$ be a smooth two vector on $M$, then the cotangent
  bundle
  is provided with the almost Lie
  algebroid structure $(T^{\star}M,\sharp,[.,.]_{\Pi}),$ where $\sharp$ is the
  vector
  bundle
  morphism from $T^{\star}M$ to $TM$ defined by $\Pi.$ The bracket of two
   differential forms $\alpha, \beta \in \Gamma(T^{\star}M)$ is defined by

   $$[\alpha,\beta]_{\Pi} = L_{\sharp\alpha}\beta -
   L_{\sharp\beta}\alpha - d\Pi (\alpha,\beta). $$

   The Jacobi anomaly of the almost Lie algebroid above is related to the
   Schouten
   square of $\Pi.$ Examples of such structures are twisted Poisson structures.
    \\

   To conclude the author apologizes  for limiting the references above to
   those
  he has needed to prepare the present work.\\

{\bf Acknowledgements}. The author would like to thank Augustin
Banyaga, Jean-Paul Dufour and Joshua Leslie for useful
discussions. He also thanks Nguen Tien Zung for carefully reading
preliminary versions of the present work.


\bibliographystyle{amsplain}

\end{document}